



%
\magnification1200
\pretolerance=100
\tolerance=200
\hbadness=1000
\vbadness=1000
\linepenalty=10
\hyphenpenalty=50
\exhyphenpenalty=50
\binoppenalty=700
\relpenalty=500
\clubpenalty=5000
\widowpenalty=5000
\displaywidowpenalty=50
\brokenpenalty=100
\predisplaypenalty=7000
\postdisplaypenalty=0
\interlinepenalty=10
\doublehyphendemerits=10000
\finalhyphendemerits=10000
\adjdemerits=160000
\uchyph=1
\delimiterfactor=901
\hfuzz=0.1pt
\vfuzz=0.1pt
\overfullrule=5pt
\hsize=146 true mm
\vsize=8.9 true in
\maxdepth=4pt
\delimitershortfall=.5pt
\nulldelimiterspace=1.2pt
\scriptspace=.5pt
\normallineskiplimit=.5pt
\mathsurround=0pt
\parindent=20pt
\catcode`\_=11
\catcode`\_=8
\normalbaselineskip=12pt
\normallineskip=1pt plus .5 pt minus .5 pt
\parskip=6pt plus 3pt minus 3pt
\abovedisplayskip = 12pt plus 5pt minus 5pt
\abovedisplayshortskip = 1pt plus 4pt
\belowdisplayskip = 12pt plus 5pt minus 5pt
\belowdisplayshortskip = 7pt plus 5pt
\normalbaselines
\smallskipamount=\parskip
 \medskipamount=2\parskip
 \bigskipamount=3\parskip
\jot=3pt
%
%
\def\ref#1{\par\noindent\hangindent2\parindent
 \hbox to 2\parindent{#1\hfil}\ignorespaces}
%
%
\font\tenss=cmss10
\font\sevenss=cmss8 at 7pt
\font\fivess=cmss8 at 5pt
\newfam\ssfam %
\textfont\ssfam=\tenss
\scriptfont\ssfam=\sevenss
\scriptscriptfont\ssfam=\fivess
%
%
%
%
%
%
%
%
%
\catcode`\_=11
\def\suf_fix{}
\def\scaled_rm_box#1{%
 \relax
 \ifmmode
   \mathchoice
    {\hbox{\tenrm #1}}%
    {\hbox{\tenrm #1}}%
    {\hbox{\sevenrm #1}}%
    {\hbox{\fiverm #1}}%
 \else
  \hbox{\tenrm #1}%
 \fi}
\def\suf_fix_def#1#2{\expandafter\def\csname#1\suf_fix\endcsname{#2}}
\def\I_Buchstabe#1#2#3{%
 \suf_fix_def{#1}{\scaled_rm_box{I\hskip-0.#2#3em #1}}
}
\def\rule_Buchstabe#1#2#3#4{%
 \suf_fix_def{#1}{%
  \scaled_rm_box{%
   \hbox{%
    #1%
    \hskip-0.#2em%
    \lower-0.#3ex\hbox{\vrule height1.#4ex width0.07em }%
   }%
   \hskip0.50em%
  }%
 }%
}
\I_Buchstabe B22
\rule_Buchstabe C51{34}
\I_Buchstabe D22
\I_Buchstabe E22
\I_Buchstabe F22
\rule_Buchstabe G{525}{081}4
\I_Buchstabe H22
\I_Buchstabe I20
\I_Buchstabe K22
\I_Buchstabe L20
\I_Buchstabe M{20em }{I\hskip-0.35}
\I_Buchstabe N{20em }{I\hskip-0.35}
\rule_Buchstabe O{525}{095}{45}
\I_Buchstabe P20
\rule_Buchstabe Q{525}{097}{47}
\I_Buchstabe R21 
\rule_Buchstabe U{45}{02}{54}
\suf_fix_def{Z}{\scaled_rm_box{Z\hskip-0.38em Z}}
\catcode`\"=12
\newcount\math_char_code
\def\suf_fix_math_chars_def#1{%
 \ifcat#1A
  \expandafter\math_char_code\expandafter=\suf_fix_fam
  \multiply\math_char_code by 256
  \advance\math_char_code by `#1
  \expandafter\mathchardef\csname#1\suf_fix\endcsname=\math_char_code
  \let\next=\suf_fix_math_chars_def
 \else
  \let\next=\relax
 \fi
 \next}
%
%
%
%
\def\font_fam_suf_fix#1#2 #3 {%
 \def\suf_fix{#2}
 \def\suf_fix_fam{#1}
 \suf_fix_math_chars_def #3.
}
\font_fam_suf_fix
 0rm
 ABCDEFGHIJKLMNOPQRSTUVWXYZabcdefghijklmnopqrstuvwxyz
\font_fam_suf_fix
 2scr
 ABCDEFGHIJKLMNOPQRSTUVWXYZ
\font_fam_suf_fix
 \slfam sl
 ABCDEFGHIJKLMNOPQRSTUVWXYZabcdefghijklmnopqrstuvwxyz
\font_fam_suf_fix
 \bffam bf
 ABCDEFGHIJKLMNOPQRSTUVWXYZabcdefghijklmnopqrstuvwxyz
\font_fam_suf_fix
 \ttfam tt
 ABCDEFGHIJKLMNOPQRSTUVWXYZabcdefghijklmnopqrstuvwxyz
\font_fam_suf_fix
 \ssfam
 ss
 ABCDEFGHIJKLMNOPQRSTUVWXYZabcdefgijklmnopqrstuwxyz
\catcode`\_=8
\def\Fdss{{\fam\ssfam I\mkern -2.5mu F}}%
\def\Ndss{{\fam\ssfam I\mkern -2.5mu N}}%
\def\Zdss{{\fam\ssfam Z\mkern-8.1mu Z}}%
%
%
%
%
\font\teneuf=eufm10 
\font\seveneuf=eufm7
\font\fiveeuf=eufm5
\newfam\euffam \def\euf{\fam\euffam\teneuf} 
\textfont\euffam=\teneuf \scriptfont\euffam=\seveneuf
\scriptscriptfont\euffam=\fiveeuf

\def\Mfr{{\euf M}}

\def\Pfr{{\euf P}}

\input xypic.tex

\def\Hom{\mathop{\rm Hom}\nolimits}
\def\Ext{\mathop{\rm Ext}\nolimits}
\def\E{\mathop{\Escr xt}\nolimits}

\parindent=0pt

\centerline{\bf On the dimension theory of skew power series rings}

\centerline{by Peter Schneider and Otmar Venjakob}

\bigskip

{\bf Introduction}

The purpose of this paper is twofold. First of all we will set up
a general notion of skew power series rings. In a skew power
series ring in one variable $t$ over a ring $R$ of coefficients
the way the variable $t$ commutes with a coefficient $a \in R$
should be subject to a relation of the form
$$
ta = \sigma(a)t + \delta(a)
$$
where $\sigma$ is an automorphism of the ring $R$ and $\delta$ is
a certain kind of derivation on $R$. This is modelled on the well
known notion of skew polynomial rings (cf.\ [MCR] 1.2). But it is
clear that in the context of power series a convergence issue will
arise. It therefore seems quite natural to work in the context of
topological rings. Having our applications in mind we will place
ourselves in the context of pseudocompact rings. In the first
section we will see that to have a well defined ring structure on
skew power series subjected to the above relation essentially
comes down to a certain nilpotency condition on the derivation
$\delta$. We also will construct a natural filtration on the skew
power series ring such that the associated graded ring is a skew
polynomial ring (in fact with zero derivation). This provides a
basic tool to establish ring theoretic properties of skew power
series rings and extends work of the second author in [Ve2].

Suppose now that $S$ is a noetherian skew power series ring over
the noetherian ring $R$. Our second theme in sections 2 and 3 is
the investigation of the class of $S$-modules which are finitely
generated as $R$-modules. In Prop.\ 2.2 we establish a short exact
sequence for any such module $M$ which presents $M$ as the
cokernel of a twisted endomorphism of the $S$-module $S \otimes_R
M$. This fact will explored in two ways. First, through a careful
analysis of the associated long exact Ext-sequence we will show in
Prop.\ 3.1 that the groups $\Ext_S^\ast(M,S)$ and $\Ext_R^{\ast
-1}(M,R)$ are (up to a twist by $\sigma$) naturally isomorphic. In
the case where the rings $R$ and $S$ are Auslander regular this
means that the geometric intuition that the codimension of $M$ as
an $S$-module is one higher than the codimension of $M$ as an
$R$-module is correct. Second, under an additional condition on
the form of the automorphism $\sigma$ it allows to show that the
class of $M$ in the Grothendieck group of modules of codimension
at most one less vanishes. In fact, this can be formulated (Prop.\
3.4) as the vanishing of a certain natural map on higher
$K$-theory, which is in the spirit of the Gersten conjecture for
commutative regular local rings (cf.\ [Ger]).

In section 4 we show that all these results apply to Iwasawa
algebras of $p$-adic Lie groups. The investigation of these rings,
which play an increasingly important role in number theory,
certainly was our original motivation. In the final section we
briefly discuss how our results extend to skew Laurent series
rings.

We remark that everything that follows holds true in an analogous
way, and with much simpler proofs, for skew polynomial rings. But
whereas the analog of our short exact sequence in Prop.\ 2.2 is
contained in [MCR] 7.5.2 its application to the study of
Ext-groups in section 3 seems to have gone unnoticed in the
literature. So the analog of our Prop.\ 3.1 constitutes a
generalization of Rees' lemma in [Ree]. In fact, we later will use
a result from [LVO] (their Thm.\ III.3.4.6) which is based on such
a generalization but which is stated in loc.\ cit.\ with an
incomplete proof.

The first author acknowledges the support of the University of
Chicago where part of this paper was written.

\medskip

{\bf 0. Notations and reminders}

\medskip

For any (unital) ring $R$ we let $\Mfr(R)$ denote the abelian category
of left (unital) $R$-modules. If $R$ is left noetherian then the
finitely generated left $R$-modules form a full abelian subcategory
$\Mfr_{fg}(R)$ of $\Mfr(R)$; we let $G(R)$ be the Grothendieck group
of the category $\Mfr_{fg}(R)$ with respect to short exact sequences.
By ${\rm Jac}(R)$ we denote the Jacobson radical of $R$.

A left pseudocompact ring $R$ is a (unital) topological ring which
is Hausdorff and complete and which has a fundamental system
$(L_i)_{i \in I}$ of open neighbourhoods of zero consisting of
left ideals $L_i \subseteq R$ such that $R/L_i$ is of finite
length as an $R$-module. A pseudocompact left $R$-module is a
(unital) topological $R$-module which is Hausdorff and complete
and which has a fundamental system $(M_i)_{i \in I}$ of open
neighbourhoods of zero consisting of submodules $M_i \subseteq M$
such that the $R$-module $M/M_i$ is of finite length. We let
$\Pfr\Mfr(R)$ denote the category of pseudocompact left
$R$-modules with continuous $R$-linear maps. The category
$\Pfr\Mfr(R)$ is abelian with exact projective limits and the
forgetful functor $\Pfr\Mfr(R) \longrightarrow \Mfr(R)$ is
faithful and exact and commutes with projective limits ([Ga1] IV.3
Thm.\ 3, [vGB] Prop.\ 3.3). This implies easily that any finitely
generated submodule of a pseudocompact module is closed and is
pseudocompact in the subspace topology. An arbitrary direct
product of pseudocompact modules $R$-modules is pseudocompact. A
pseudocompact $R$-module $M$ will be called topologically free if
it is topologically isomorphic to a pseudocompact module of the
form $\prod_{i \in I} R$; the set of elements of $M$ corresponding
to the elements $(\ldots,0,1,0,\ldots)$, with the $1$ in the i-th
place, then is called a pseudobasis of $M$. Topologically free
pseudocompact modules are projective objects in the category
$\Pfr\Mfr(R)$; one easily concludes that $\Pfr\Mfr(R)$ has enough
projective objects ([Bru] Cor. 1.3 and Lemma 1.6). In particular,
for any two pseudocompact $R$-modules $M$ and $N$, Ext-functors
$\E^\ast_R(M,N)$ may be defined by left deriving the functor
$\Hom_{\Pfr\Mfr(R)}(.,N)$. Suppose that $R$ is left noetherian.
Then any finitely generated left $R$-module carries a unique
pseudocompact topology. Hence we have a natural fully faithful and
exact embedding $\Mfr_{fg}(R) \longrightarrow \Pfr\Mfr(R)$. Given
a finitely generated $R$-module $M$ and an arbitrary pseudocompact
$R$-module $N$ one easily shows, by using a resolution of $M$ by
finitely generated projective $R$-modules, that one has
$$
\E^\ast_R(M,N) = \Ext^\ast_R(M,N) \ ,
$$
where the right hand side denotes the usual Ext-groups of any two
modules in $\Mfr(R)$.

There are obvious ``right'' versions of everything above. If we are
dealing with these we usually will indicate this by adding a
superscript ``r'' to the notation; e.g., $\Mfr^r(R)$ denotes the
category of right $R$-modules. By an ideal we always mean a 2-sided
ideal, and a ring will be called noetherian if it is left and right
noetherian.

\smallskip

{\it Remark 0.1:} Suppose that $R$ is left noetherian; we then have:

i. If $R$ is left pseudocompact then\hfill\break
 -- the ring $R/{\rm Jac}(R)$ is left artinian, and\hfill\break
 -- the pseudocompact topology on $R$ is the ${\rm Jac}(R)$-adic topology;

ii. if there is an ideal $I \subseteq R$ such that $R/I$ is left
artinian and $R$ is $I$-adically complete then $R$ is left
pseudocompact in the $I$-adic topology.

Proof: i. See [vGB] Cor. 3.14. ii. Each $I^n/I^{n+1}$ is a finitely
generated left module over the left artinian ring $R/I$ and hence is
of finite length. It follows inductively that each $R/I^n$ is a left
$R$-module of finite length.

\smallskip

{\it Remark 0.2:} Suppose the left pseudocompact ring $R$ is
noetherian; then $R$ also is right pseudocompact (in the same
topology).

Proof: This follows from Remark 0.1 since the noetherian and left
artinian ring $R/I$ necessarily is right artinian as well ([CH]
p.55).

\smallskip

Suppose that $R \longrightarrow S$ is a continuous (unital)
homomorphism of noetherian pseudocompact rings. For any pseudocompact
left $R$-module $M$ we define its (completed) {\it base change} to $S$
by
$$
S\, \widehat{\otimes}_R\, M := \mathop{\rm
lim}\limits_{\mathop{\longleftarrow}\limits_{i,L}}\, S/{\rm Jac}(S)^i
\otimes_R M/L
$$
where $i$ runs over $\Ndss$ and $L$ runs over all open submodules of
$M$. Since $M/L$ is of finite length over $R$ we find an integer $j
\geq 0$ such that ${\rm Jac}(R)^jM \subseteq L$. By possibly
increasing $j$ we have
$$
S/{\rm Jac}(S)^i \otimes_R M/L = S/{\rm Jac}(S)^i
\otimes_{R/{\rm Jac}(R)^j} M/L\ .
$$
Because $M/L$ is finitely generated over the artinian ring $R/{\rm
Jac}(R)^j$ it follows that $S/{\rm Jac}(S)^i \otimes_{R/{\rm
Jac}(R)^j} M/L$ is finitely generated over the artinian ring
$S/{\rm Jac}(S)^i$ and hence is of finite length over $S$.
Therefore $S\, \widehat{\otimes}_R\, M$ is a projective limit of
finite length $S$-modules and as such is pseudocompact. We have
constructed in this way a functor
$$
S\, \widehat{\otimes}_R\, . : \Pfr\Mfr(R) \longrightarrow \Pfr\Mfr(S)
\ .
$$
It is right exact since projective limits in $\Pfr\Mfr(S)$ are exact.
If $M$ is finitely generated then, by looking at a finite presentation
of $M$ over $R$, it is rather obvious that one has
$$
S\, \widehat{\otimes}_R\, M = S \otimes_R M \ .
$$

\smallskip

{\it Remark 0.3:} If $S$ is pseudocompact and topologically free as a
right $R$-module then the functor $S\, \widehat{\otimes}_R\, .$ is
exact.

Proof: As a topological abelian group $S\, \widehat{\otimes}_R\,
M$ only depends on the structure of $S$ as a pseudocompact right
$R$-module. Under our additional assumption we may write $S$ as
the projective limit $S = \mathop{\rm lim}\limits_{\longleftarrow}
S_i$ of right factor $R$-modules $S_i$ of $S$ which are finitely
generated and free over $R$. By the argument in the proof of [Bru]
Lemma A.4 the natural map
$$
S\, \widehat{\otimes}_R\, M \longrightarrow
\mathop{\rm lim}\limits_{\longleftarrow}\, (S_i \otimes_R M)
$$
is bijective. Since $S_i$ is free over $R$ and since the transition
maps in the projective system $\{S_i \otimes_R M\}_i$ are surjective
(for any $M$) the right hand side of the above bijection is an exact
functor in $M$.

\medskip

{\bf 1. Skew power series rings}

\medskip

In order to motivate our later construction of skew power series
rings we first recall briefly the well known notion of a skew
polynomial ring (cf.\ [MCR] 1.2). Let $R$ be an arbitrary unital
ring. We suppose given an automorphism $\sigma$ of the ring $R$ as
well as a left $\sigma$-derivation $\delta : R \rightarrow R$
which is an additive map satisfying
$$
\delta(ab) = \delta(a)b + \sigma(a)\delta(b) \qquad\qquad\hbox{for
any}\ a,b \in R\ .
$$
The left $R$-module of all (left) polynomials $a_0 + a_1t + \ldots
+ a_mt^m$ over $R$ in the variable $t$ carries a unique unital
ring structure which extends the $R$-module structure and
satisfies
$$
ta = \sigma(a)t + \delta(a) \qquad\qquad\hbox{for any}\ a \in R\ .
$$
This ring is denoted by $R[t;\sigma,\delta]$. Put $\sigma' :=
\sigma^{-1}$ and $\delta' := - \delta\sigma^{-1}$. Then $\delta'$
is a right $\sigma'$-derivation, i.e., an additive map satisfying
$$
\delta'(ab) = \delta'(a)\sigma'(b) + a\delta'(b)
\qquad\qquad\hbox{for any}\ a,b \in R\ .
$$
The skew polynomial ring $R[t;\sigma,\delta]$ can alternatively be
described as the right $R$-module of all (right) polynomials $a_0
+ ta_1 + \ldots + t^ma_m$ over $R$ with the multiplication
determined by the relation
$$
at = t\sigma'(a) + \delta'(a) \qquad\qquad\hbox{for any}\ a \in R\
.
$$
To make this explicit we let $M_{k,l}(Y,Z)$, for any integers $k,l
\geq 0$, denote the sum of all noncommutative monomials in two
variables $Y,Z$ with $k$ factors $Y$ and $l$ factors $Z$. In
$R[t;\sigma,\delta]$ one then has the formulas
$$
\sum_{i \geq 0} t^ia_i = \sum_{j \geq 0} \Big(\sum_{i \geq j}
M_{i-j,j}(\delta,\sigma)(a_i)\Big) t^j \leqno{(1)}
$$
and
$$
\sum_{j \geq 0} a_jt^j = \sum_{i \geq 0} t^i \Big(\sum_{j \geq i}
M_{j-i,i}(\delta',\sigma')(a_j)\Big)\ .\leqno{(2)}
$$
Furthermore, the multiplication is explicitly given by
$$
(\sum_{i \geq 0} t^ia_i)(\sum_{k \geq 0} t^ka_k) = \sum_{m \geq 0}
t^m \Big(\sum_{n=0}^m \sum_{k \geq n}
M_{k-n,n}(\delta',\sigma')(a_{m-n})b_k\Big) \leqno{(3)}
$$
and
$$
(\sum_{j \geq 0} a_jt^j)(\sum_{l \geq 0} b_lt^l) = \sum_{m \geq 0}
\Big(\sum_{n=0}^m \sum_{j \geq n}
a_jM_{j-n,n}(\delta,\sigma)(b_{m-n})\Big) t^m \ , \leqno{(4)}
$$
respectively.

After this reminder we always in this section let $R$ be a
noetherian pseudocompact ring, $\sigma$ a topological automorphism
of $R$, and $\delta : R \rightarrow R$ a continuous left
$\sigma$-derivation. We define $S$ to be the left $R$-module of
all (left) formal power series $\sum_{i \geq 0} a_it^i$ over $R$
in the variable $t$ and view it as a pseudocompact $R$-module with
respect to the obvious direct product topology. Then a formal
power series $ x = \sum_{i \geq 0} a_it^i$ can be considered as an
expansion of the element $x$ into a convergent sum. The skew
polynomial ring $R[t;\sigma,\delta]$ is a dense submodule of $S$.
To extends its ring structure by continuity to $S$ we need to
ensure that the sums which form the coefficients on the right hand
sides of the formulas (1) - (4) and which then are infinite
converge. We let $\Ascr_k$ denote the set of all noncommutative
monomials in three variables $Y$,$Z$, and $Z'$ with exactly $k$
factors $Y$, and we put $\Ascr_{\geq l} := \bigcup_{k \geq l}
\Ascr_k$.

\smallskip

{\bf Definition:} {\it A left or right $\sigma$-derivation
$\widetilde{\delta} : R \rightarrow R$ is called
$\sigma$-nilpotent if for any $n \geq 1$ there is an $m \geq 1$
such that
$$
M(\widetilde{\delta},\sigma,\sigma^{-1})(R) \subseteq {\rm
Jac}(R)^n \qquad\hbox{for any}\ M \in \Ascr_{\geq m}\ .
$$
}

\smallskip

{\it Remark 1.1:} i. If $\widetilde{\delta}$ and $\sigma$ commute
then $M(\widetilde{\delta},\sigma,\sigma^{-1})(R) =
\widetilde{\delta}^k(R)$ for $M \in \Ascr_k$ and hence
$\widetilde{\delta}$ is $\sigma$-nilpotent if and only if it is
topologically nilpotent.

ii. If $\widetilde{\delta}$ satisfies $\widetilde{\delta}(R)
\subseteq {\rm Jac}(R)$ and $\widetilde{\delta}({\rm Jac}(R))
\subseteq {\rm Jac}(R)^2$ then $\widetilde{\delta}({\rm Jac}(R)^n)
\subseteq {\rm Jac}(R)^{n+1}$ holds true for any $n \geq 0$ and,
because of $\sigma^{\pm 1}({\rm Jac}(R)) = {\rm Jac}(R)$, it
follows that $\widetilde{\delta}$ is $\sigma$-nilpotent.

\smallskip

We assume from now on that $\delta$ is $\sigma$-nilpotent. As
before we put $\sigma' := \sigma^{-1}$ and $\delta' := -
\delta\sigma^{-1}$. Clearly $\delta'$ is $\sigma'$-nilpotent. One
easily checks that by reading (4) as a definition we obtain a
continuous multiplication on $S$ with respect to which $S$ is a
left pseudocompact ring (compare the reasoning in [Ga2] 0.4).
Since also the formulas (1) - (3) hold now more generally in $S$
we see that $S$ is right pseudocompact as well and that $\{t^i : i
\geq 0\}$ is a pseudobasis of $S$ as a left and as a right
$R$-module. We call $R[[t;\sigma,\delta]] := S$ a (formal) {\it
skew power series ring} over $R$.

\smallskip

{\it Remark 1.2:} If $\sigma$ and $\delta$ commute with each other
then the map
$$
\matrix{
\widehat{\sigma} : & \hfill S & \longrightarrow & S \hfill\cr
 & \sum_i t^ia_i & \longmapsto & \sum_i t^i\sigma(a_i) }
$$
is a topological automorphism of the pseudocompact ring $S$.

Proof: The map obviously is additive, bijective, and a
homeomorphism. It remains to check its multiplicativity on
elements of the form $x = t^ia$ and $y = t^jb$ with $a,b \in R$.
By a straightforward induction it in fact suffices to consider the
case where $j = 1$ and $b = 1$. We compute
$$
\widehat{\sigma}(t^iat) = \widehat{\sigma}(t^i(t\sigma^{-1}(a) -
\delta(\sigma^{-1}(a)))) = t^{i+1}a - t^i\delta(a) = t^i\sigma(a)t =
\widehat{\sigma}(t^ia)\widehat{\sigma}(t) \ .
$$


\smallskip

In order to transfer properties of the ring $R$ to the ring $S$ we
introduce a certain filtration on $R$. For $k \in \Ndss$ let $P_k$
denote the set of all sequences $\underline{m} = (m_1,\ldots,m_r)$
of natural numbers $m_i > 0$ and varying length $r$ such that $m_1
+ \ldots + m_r = k$. We define $\Delta_0 := \{1_R\}$ and
$$
\Delta_k := \sum_{\underline{m} \in P_k}\sum_{M_i\in \Ascr_{m_i}}
M_1(\delta,\sigma,\sigma^{-1})(R)\cdot\ldots\cdot
M_r(\delta,\sigma,\sigma^{-1})(R)
$$
for $k\geq 1.$ The following properties are easily verified by
induction:

(a)\quad $\Delta_k\cdot \Delta_l \subseteq \Delta_{k+l}$,

(b)\quad $\sigma (\Delta_k) = \Delta_k$,

(c)\quad $\delta(\Delta_k)\subseteq \Delta_{k+1},\;
\delta(R\Delta_k)\subseteq R\Delta_{k+1},\; \delta(\Delta_k
R)\subseteq \Delta_{k+1} R$,

(d)\quad $R\Delta_{k+1}\subseteq R\Delta_k$,

(e)\quad $R\Delta_k=\Delta_k R$.

In particular, the $I_k:=R\Delta_k$, for $k\geq 0$, form a
descending series of 2-sided ideals in $R$ with $I_k\cdot
I_l\subseteq I_{k+l}$. The associated graded ring is denoted, as
usual, by $gr_I R$. The automorphism $\sigma$ induces an
automorphism $\overline{\sigma}$ of $gr_I R$ whereas $\delta$
induces the zero map on $gr_I R$.

\smallskip

{\it Remark 1.3:} i. If $\delta$ is $\sigma$-nilpotent with
$\delta(R) \subseteq {\rm Jac}(R)$ then we have $\bigcap I_k=0$;

ii. if $\bigcap I_k=0$ then $\delta$ is $\sigma$-nilpotent and $R
= \mathop{\rm lim}\limits_{\longleftarrow} R/I_k$.

Proof: i. This is obvious. ii. First of all we note that
$M(\delta,\sigma,\sigma^{-1})(R) \subseteq \Delta_k \subseteq I_k$
for $M \in \Ascr_k$. On the other hand, since $R$ is noetherian
the ideals $I_k$ are closed in $R$. It therefore follows from
[Ga1] IV.3. Prop.\ 11 that for any $n \geq 1$ there is a $k \geq
1$ such that $I_k \subseteq {\rm Jac}(R)^n$. Moreover, [Ga1] IV.3
Prop.\ 10 says that the canonical map $R \longrightarrow
\mathop{\rm lim}\limits_{\longleftarrow} R/I_k$ is an isomorphism
of pseudocompact rings.

\smallskip

The filtration $I_k$ of $R$ (with $I_k := R$ for $k < 0$) induces
a filtration $J_k$ by
$$
J_k:=\prod_{i \geq 0} I_{k-i}t^i \ .
$$
Indeed we have the following.

\smallskip

{\bf Lemma 1.4:} {\it i. Each $J_k$ is a closed 2-sided ideal in
$S$;

ii. $J_k\cdot J_l\subseteq J_{k+l},$;

iii. if $\;\bigcap I_k = 0$ then $S = \mathop{\rm
lim}\limits_{\longleftarrow} S/J_k$;

iv. for the associated graded ring $gr_J S$ we have $gr_J S \cong
(gr_IR)[\overline{t};\overline{\sigma}]$ where $\overline{t}$
denotes the principal symbol of $t$ in $gr_J S$.}

Proof: i. and ii. are direct consequences of the above properties
of $I_k$ and the earlier formula (4). iii. If $\bigcap I_k = 0$
then also $\bigcap J_k = 0$ and we may apply again [Ga1] IV.3
Prop.\ 10.iv. This is straightforward.

\smallskip

{\bf Lemma 1.5:} {\it Suppose that $\delta$ is $\sigma$-nilpotent
with $\delta(R) \subseteq {\rm Jac}(R)$ and that $gr_I R$ is
noetherian, resp.\ Auslander regular; then $S$ is noetherian,
resp.\ Auslander regular.}

Proof: (The notion of Auslander regularity will be recalled in
section 3.) By Remark 1.3.i and Lemma 1.4.iii the ring $S$ is
complete with respect to the filtration $J_k$. Moreover, it is
well known that the skew polynomial ring $gr_J S \cong gr_I
R[\overline{t};\overline{\sigma}]$ over the noetherian, resp.\
Auslander regular, ring $gr_I R$ is noetherian, resp.\ Auslander
regular, as well (cf.\ [MCR] 1.2.9 and [LVO] III.3.4.6). Finally,
it is a general fact that a complete filtered ring is noetherian,
resp.\ Auslander regular, if its associated graded ring has this
property (cf.\ [LVO] I.1.2.3 and III.2.2.5).

\smallskip

Since in our main application in section 4 the maps $\sigma$ and
$\delta$ commute with each other we want to mention an example
where this is not the case. Let $R = R_0[[X]]$ be the commutative
formal power series ring in one variable $X$ over the ring $R_0 =
\Zdss_p$ or $R_0 = \Fdss_p$. Its Jacobson radical ${\rm Jac}(R)$
is the ideal generated by $p$ and $X$. We fix a unit $u \in
R^\times$ and let $\sigma$ denote the unique continuous
automorphism of $R$ such that $\sigma(X) = uX$ and $\sigma|R_0 =
id$. We also fix an element $F \in {\rm Jac}(R)^2$. There is a
unique continuous $\sigma$-derivation $\delta$ on $R$ such that
$\delta(X) = F$ and $\delta|R_0 = 0$. It satisfies $\delta({\rm
Jac}(R)^j) \subseteq {\rm Jac}(R)^{j+1}$ for any $j \geq 0$ and
therefore is $\sigma$-nilpotent. Hence the corresponding skew
power series ring $R[[t;\sigma,\delta]]$ is defined. Obviously,
$\sigma$ and $\delta$ do not commute with each other in general.
To be completely explicit let $R_0 = \Fdss_p$, $F := X^p$, and $u
\in \Fdss_p\setminus\{0,1\}$. Then $I_k = RX^{pk}$; hence $gr_I R$
is the polynomial ring in one variable over
$\Fdss_p[X]/\Fdss_p[X]X^p$.

\smallskip

In our main application later on the situation at first is a
little bit different. The data $(R,\sigma,\delta)$ are as before
but $\delta$ is not assumed to be $\sigma$-nilpotent. Instead $R$
is contained in another noetherian pseudocompact ring $S'$ such
that:

-- $S'$ is pseudocompact and topologically free as a left as well
as a right $R$-module;

-- there is an element $t \in S'$ such that the powers $\{t^i : i
\geq 0\}$ form a pseudobasis of $S'$ as a left and as a right
$R$-module and such that
$$
ta = \sigma(a)t + \delta(a) \qquad\qquad\hbox{for any}\ a \in R\ .
$$

\smallskip

{\bf Lemma 1.6:} {\it Suppose that we are in the above situation
(in particular, $S'$ is noetherian) and that $\sigma$ and $\delta$
commute with each other; we then have:

i. $\delta^n(a) = \sum_{i=0}^n {n\choose i} (-1)^i
t^{n-i}\sigma^i(a)t^i$ for any $a \in R$ and any integer $n \geq
0$;

ii. for any open left ideal $L \subseteq R$ there is an $n \geq 0$
such that $\delta^j(R) \subseteq L$ for any $j \geq n$.}

Proof: i. This is proved by an easy induction starting from the
defining equation for a $\sigma$-derivation.

ii. The topology of $R$ is induced by the topology of $S'$. Since
$S'$ is noetherian we therefore find a $k \geq 0$ such that ${\rm
Jac}(S')^k \cap R \subseteq L$. Since ${\rm Jac}(S')^k$ is open in
$S'$ all but finitely many members of the pseudobasis $\{t^i : i
\geq 0\}$ must lie in ${\rm Jac}(S')^k$. Hence there is an $m \geq
0$ such that $S't^mS' \cap R \subseteq {\rm Jac}(S')^k \cap R
\subseteq L$. But it is an immediate consequence of the formula in
i. that $\delta^{2m}(R) \subseteq S't^mS' \cap R$.


\smallskip

We see that in the situation of the lemma $\delta$ necessarily is
topologically nilpotent, hence $\sigma$-nilpotent by Remark 1.1.i,
and that therefore $S'$ is isomorphic to the skew power series
ring $R[[t;\sigma,\delta]]$.

\medskip

{\bf 2. A short exact sequence}

\medskip

Throughout this paper we fix a skew power series ring $S =
R[[t;\sigma,\delta]]$ over a noetherian pseudocompact ring $R$ as
constructed in the previous section. In particular, $\delta$ is
$\sigma$-nilpotent. We also assume that $S$ is noetherian.

\smallskip

{\bf Lemma 2.1:} {\it For any module $M$ in $\Mfr_{fg}(S)$ we have
$\bigcap_{k \geq 0} t^kM = \{0\}$. }

Proof: By [vGB] Cor. 3.14 the pseudocompact topology on $M$ is the
${\rm Jac}(S)$-adic one. It follows that $\bigcap_{k \geq 0} {\rm
Jac}(S)^kM = \{0\}$. Since ${\rm Jac}(S)^k$ is open in $S$ all but
finitely many members of the pseudobasis $\{t^i : i \geq 0\}$ lie
in ${\rm Jac}(S)^k$. Hence $\bigcap_{k \geq 0} t^kM \subseteq
\bigcap_{k \geq 0} {\rm Jac}(S)^kM = \{0\}$.

\smallskip

We claim that the forgetful functor induces an inclusion of abelian
categories
$$
\Pfr\Mfr(S) \subseteq \Pfr\Mfr(R) \ .
$$
Let $M$ be a pseudocompact left $S$-module and let $M_0 \subseteq
M$ be an open submodule. Then $M/M_0$ is of finite length as
$S$-module. Hence we have ${\rm Jac}(S)^i M \subseteq M_0$ for
some $i > 0$ so that in fact $M/M_0$ is an $S/{\rm
Jac}(S)^i$-module of finite length. But by Remark 0.1.i, since $S$
is assumed to be noetherian, the $R$-module $S/{\rm Jac}(S)^i$ and
consequently the $R$-module $M/M_0$ is of finite length. It
follows that $M$ is pseudocompact as an $R$-module. Clearly, the
base extension functor $S\, \widehat{\otimes}_R\, .$ is left
adjoint to the above inclusion functor and hence preserves
projective objects.

For any left $R$-module $M$ the twisted left $R$-module $^\sigma M$
has the same underlying additive group as $M$ but with $R$ acting
through the automorphism $\sigma^{-1}$. If $M$ is pseudocompact then
$^\sigma M$ obviously is pseudocompact as well.

We also introduce the abelian category $\Mfr_R(S)$ of those left
$S$-modules which are finitely generated as $R$-modules.

\smallskip

{\bf Proposition 2.2:} {\it For any module $M$ in $\Mfr_R(S)$ the
sequence
$$
0\longrightarrow S \otimes_R {^\sigma M}
\mathop{\longrightarrow}\limits^{\kappa} S \otimes_R M
\mathop{\longrightarrow}\limits^{\mu} M \longrightarrow 0
$$
with $\kappa(x \otimes m) := xt \otimes m - x \otimes tm$ and $\mu(x
\otimes m) := xm$ is an exact sequence of finitely generated
$S$-modules.}

Proof: Obviously all $S$-modules in the sequence under consideration
are finitely generated. To see that $\kappa$ is well defined fix an $a
\in R$. Then
$$
\eqalign{
\kappa(xa \otimes m) & = xat \otimes m - xa \otimes tm\cr
 & = xat \otimes m - x \otimes atm\cr
 & = xt\sigma^{-1}(a) \otimes m - x \otimes t\sigma^{-1}(a)m\cr
 & = xt \otimes \sigma^{-1}(a)m - x \otimes t\sigma^{-1}(a)m\cr
 & = \kappa(x \otimes \sigma^{-1}(a)m)\ .\cr
}
$$
Clearly, $\kappa$ and $\mu$ are $S$-module maps such that
$\mu\circ\kappa = 0$, and $\mu$ is surjective. It remains to show that
$\kappa$ is injective and that the kernel of $\mu$ is contained in the
image of $\kappa$.

It is immediate that the kernel of $\mu$ is generated additively by
the elements of the form
$$
x \otimes m - 1 \otimes xm \qquad\hbox{for}\ x \in S\
\hbox{and}\ m \in M\ .
$$
As a finitely generated module $S \otimes_R M$ carries a natural
pseudocompact topology with respect to which any submodule, being
finitely generated as well, is closed. It therefore follows from the
power series expansion of $x$ that ${\rm ker}(\mu)$ as a $S$-module is
generated by the elements
$$
t^i \otimes m - 1 \otimes t^im \qquad\hbox{for}\ i \geq 0\ \hbox{and}\
m \in M\ .
$$
But
$$
t^i \otimes m - 1 \otimes t^im = \sum_{j=0}^{i-1} t^j(t \otimes
t^{i-j-1}m - 1 \otimes tt^{i-j-1}m)\ .
$$
Hence in fact the elements $t \otimes m - 1 \otimes tm = \kappa(1
\otimes m) \in {\rm im}(\kappa)$ generate ${\rm ker}(\mu)$.

Finally, to establish the injectivity of $\kappa$ it is convenient to
make the identifications
$$
\matrix{
\hfill S \otimes_R {^\sigma M} &
\mathop{\longrightarrow}\limits^\simeq & \prod_{i \geq 0} M \hfill
& \hbox{and} & \hfill S \otimes_R M &
\mathop{\longrightarrow}\limits^\simeq & \prod_{i \geq 0} M
\hfill\cr\cr
\hfill (\sum_{i \geq 0} t^ia_i) \otimes m & \longmapsto & (\sigma^{-1}(a_i)m)_i
& & (\sum_{i \geq 0} t^ia_i) \otimes m & \longmapsto & (a_im)_i
\hfill }
$$
which are possible since $M$ is finitely generated over $R$. A
straightforward computation shows that the map $\kappa$ under this
identification is given by $(m_i)_i \longmapsto (m_{i-1} -
tm_i)_i$ (where $m_{-1} := 0$). An element in the kernel of
$\kappa$ therefore corresponds to a tuple $(m_i)_i$ such that $m_i
= tm_{i+1}$ for any $i \geq -1$. Hence all $m_i$ lie in the
intersection $\bigcap_{k \geq 0} t^kM$ which by Lemma 2.1 is equal
to zero.

\smallskip

{\it Remark 2.3:} i. Suppose that $\sigma$ is of the form
$\sigma(.) = \gamma . \gamma^{-1}$ for some unit $\gamma \in S$.
Then  for any $M$ in $\Mfr_R(S)$ the map $^\sigma M
\mathop{\longrightarrow}\limits^{\gamma} M$ is an isomorphism of
$R$-modules. Prop.\ 2.2 therefore implies that we have a short
exact sequence of finitely generated $S$-modules of the form
$$
0 \longrightarrow S \otimes_R M \longrightarrow S \otimes_R M
\longrightarrow M \longrightarrow 0 \ .
$$
In particular, the class of $M$ in the Grothendieck group $G(S)$ is
zero.

ii. For any module $M$ in $\Mfr^r_R(S)$ there is a corresponding exact
sequence
$$
0 \longrightarrow M^{\sigma^{-1}} \otimes_R S
\mathop{\longrightarrow}\limits^{\kappa^r}
M \otimes_R S
\mathop{\longrightarrow}\limits^{\mu^r}
M \longrightarrow 0
$$
with $\kappa^r(m \otimes x) := m \otimes tx - mt \otimes x$ and
$\mu^r(m \otimes x) := mx$.

\medskip

\goodbreak

{\bf 3. Comparing dimensions}

\medskip

In this section we want to compare the Ext-groups
$\Ext^\ast_R(M,R)$ and $\Ext^\ast_S(M,S)$ for any module $M$ in
$\Mfr_R(S)$. Since the rings $R$ and $S$ are noetherian these
Ext-groups certainly are finitely generated right $R$- and
$S$-modules, respectively. We will show that the right $R$-module
structure on $\Ext^\ast_R(M,R)^\sigma$ can be extended in a
natural way to a right $S$-module structure. In fact we start more
generally with an arbitrary finitely generated left $S$-module $M$
(which therefore is pseudocompact as $S$- and as $R$-module) and
construct first a natural right $S$-module structure on
$\Hom_{\Pfr\Mfr(R)}(M,R)^\sigma$. The latter is a right $R$-module
via
$$
f^a(m) := f(m)\sigma^{-1}(a)
$$
for any $f \in \Hom_{\Pfr\Mfr(R)}(M,R)$, $a \in R$, and $m \in M$. We
now define
$$
f^t(m) := \sigma^{-1}(f(tm) - \delta(f(m)))\ .
$$
Obviously, $f^t : M \longrightarrow R$ is additive and continuous. the
computation
$$
\eqalign{
f^t(bm) & = \sigma^{-1}(f(tbm) - \delta(f(bm)))\cr
 & = \sigma^{-1}(f(\sigma(b)tm + \delta(b)m) - \delta(bf(m)))\cr
 & = \sigma^{-1}(\sigma(b)f(tm) + \delta(b)f(m) - \delta(b)f(m) -
 \sigma(b)\delta(f(m)))\cr
 & = b\sigma^{-1}(f(tm) - \delta(f(m))) }
$$
for any $b \in R$ shows that $f^t$ also is $R$-linear. Hence $f^t \in
\Hom_{\Pfr\Mfr(R)}(M,R)$. We next check that
$$
\matrix{
[(f^{\sigma(a)})^t + f^{\delta(a)}](m) \hfill\cr
 \qquad\qquad  = \sigma^{-1}(f(tm)a - \delta(f(m)a)) +
     f(m)\sigma^{-1}(\delta(a))\hfill\cr
 \qquad\qquad = \sigma^{-1}(f(tm)a - \delta(f(m))a - \sigma(f(m))\delta(a))
     + f(m)\sigma^{-1}(\delta(a))\hfill\cr
 \qquad\qquad = \sigma^{-1}(f(tm) - \delta(f(m)))\sigma^{-1}(a)\hfill\cr
 \qquad\qquad = [(f^t)^a](m)\ .\hfill }
$$
This shows that, if we let $S_0 \subseteq S$ denote the {\it skew
polynomial} subring of all finite expressions $\sum_{i=0}^n t^ia_i$,
then by the above definitions $S_0$ acts from the right on
$\Hom_{\Pfr\Mfr(R)}(M,R)^\sigma$. Moreover, this action is functorial
with respect to maps in $\Pfr\Mfr(S)$. To see that it extends, by
continuity, to an $S$-action we first establish the following
assertion: For any $j \geq 0$ there is a $k(j) \geq 0$ such that
$$
f^{t^k}(M) \subseteq {\rm Jac}(R)^j \qquad\hbox{for any}\ k \geq k(j).
$$
For $k \geq 0$ and $0 \leq i \leq k$ we define inductively
explicit noncommutative polynomials $B_{i,k}(Y,Z,Z')$ by
$$
B_{0,k} := 1 ,\; B_{k+1,k+1} := Z^kYZ'^kB_{k,k}
$$
and
$$
B_{i,k} := B_{i,k-1} + Z^{k-1}YZ'^{k-1}B_{i-1,k-1}
\qquad\qquad\hbox{for}\ 0 < i < k \ .
$$
We note that each $B_{i,k}$ is a sum of monomials in $\Ascr_i$.
An easy induction then establishes the explicit formula
$$
(f^{t^k})(m) = \sigma^{-k}(\sum_{i=0}^k (-1)^i
B_{i,k}(\delta,\sigma,\sigma^{-1})(f(t^{k-i}m)))\leqno{(+)}
$$
for any $k \geq 0$. By the $\sigma$-nilpotence of $\delta$ there
is an $r \geq 1$ such that
$$
B(\delta,\sigma,\sigma^{-1})(R) \subseteq {\rm Jac}(R)^j
\leqno{(a)}
$$
for any finite sum $B$ of monomials in $\Ascr_{\geq r}$. Since
$\delta$ is continuous we find a sequence of integers $\ell_r
\geq\ldots\geq \ell_2 \geq \ell_1 := j$ such that
$$
\delta^0({\rm Jac}(R)^{\ell_\rho}),\ldots,\delta^{r - 1}({\rm
Jac}(R)^{\ell_\rho}) \subseteq {\rm Jac}(R)^{\ell_{\rho -1}}
\qquad\hbox{for any}\ 2 \leq \rho \leq r\ .\leqno{(b)}
$$
Because $f$ is continuous and $M$ is finitely generated there is,
by Remark 0.1.i, an $n' \geq 0$ such that
$$
f({\rm Jac}(S)^{n'}M) \subseteq {\rm Jac}(R)^{\ell_r}\ .
$$
From the proof of Lemma 2.1 we know that
$$
t^n \subseteq {\rm Jac}(S)^{n'}
$$
for some $n \geq 0$. Hence
$$
f(t^{i'}M) \subseteq {\rm Jac}(R)^{\ell_r} \qquad\hbox{for any}\
i' \geq n\ .
$$
Together with (b) and the fact that $\sigma$ is an automorphism of
$R$ this shows that
$$
B(\delta,\sigma,\sigma^{-1})(f(t^{i'}M))  \subseteq {\rm Jac}(R)^j
\qquad\hbox{for any}\ i' \geq n \leqno{(c)}
$$
whenever $B$ is a finite sum of monomials in $\Ascr_0
\cup\ldots\cup \Ascr_{r-1}$. We now choose $k \geq k(j) := r + n$.
For any index $i$ in the sum on the right hand side of the formula
$(+)$ we then have $i \geq r$ or $i < r, k - i \geq n$ which means
that the corresponding summand
$B_{i,k}(\delta,\sigma,\sigma^{-1})(f(t^{k-i}m))$ lies in ${\rm
Jac}(R)^j$ either by (a) or by (c). Hence
$$
f^{t^k}(M) \subseteq \sigma^{-k}({\rm Jac}(R)^j) = {\rm Jac}(R)^j
\ .
$$
Consider now an arbitrary element $x = \sum_{i \geq 0} t^ia_i$ in
$S$. Since $R$ is complete the above assertion implies that the
series
$$
(f^x)(m) := \sum_{i \geq 0} ((f^{t^i})^{a_i})(m)
$$
converges in $R$. It is clear that the resulting map $f^x : M
\longrightarrow R$ is $R$-linear and continuous. It is straightforward
to check that this indeed defines a right $S$-action on
$\Hom_{\Pfr\Mfr(R)}(M,R)^\sigma$ which is functorial with respect to
maps in $\Pfr\Mfr(S)$.

Let now $P. \longrightarrow M$ be a projective resolution of $M$
in $\Mfr_{fg}(S)$. Since $S$ is topologically free as a left
$R$-module this also is a projective resolution of $M$ in
$\Pfr\Mfr(R)$. It follows that the terms in the complex
$\Hom_{\Pfr\Mfr(R)}(P.,R)^\sigma$ and hence the Ext-groups
$\E^\ast_R(M,R)^\sigma$ carry a natural right $S$-module
structure.

Let us fix a module $M$ in $\Mfr_R(S)$. Then, by the above discussion,
the
$$
\Ext^j_R(M,R)^\sigma = \E^j_R(M,R)^\sigma
$$
are right $S$-modules in the category $\Mfr^r_R(S)$. Applying the
right module version of Prop.\ 2.2 to them we obtain the short
exact sequences of right $S$-modules
$$
0 \longrightarrow \Ext^j_R(M,R) \otimes_R S
\mathop{\longrightarrow}\limits^{\kappa^r}
\Ext^j_R(M,R)^\sigma \otimes_R S
\mathop{\longrightarrow}\limits^{\mu^r}
\Ext^j_R(M,R)^\sigma \longrightarrow 0 \ .
$$
On the other hand we may pass from the short exact sequence of
Prop.\ 2.2 to the associated long exact Ext-sequence and obtain
the exact sequence of right $S$-modules
$$
\matrix{
\qquad \Ext^j_S(S \otimes_R M,S) \longrightarrow
\Ext^j_S(S \otimes_R {^\sigma M},S) \longrightarrow
\Ext^{j+1}_S(M,S)\hfill\cr\cr
\qquad \longrightarrow \Ext^{j+1}_S(S \otimes_R M,S) \longrightarrow
\Ext^{j+1}_S(S \otimes_R {^\sigma M},S) \hfill }
$$
where the first and the last arrow are given by
$\Ext^\ast_S(\kappa,S)$. We claim that this latter map can naturally
be identified with the map $\kappa^r$ for the right $S$-module
$\Ext^\ast_R(M,R)^\sigma$. Hence the second exact sequence gives rise
to the short exact sequence
$$
0 \longrightarrow \Ext^j_S(S \otimes_R M,S) \longrightarrow
\Ext^j_S(S \otimes_R {^\sigma M},S) \longrightarrow
\Ext^{j+1}_S(M,S) \longrightarrow 0
$$
which moreover is naturally isomorphic to the first exact sequence
above. This establishes the following fact.

\smallskip

{\bf Proposition 3.1:} {\it For any module $M$ in $\Mfr_R(S)$ we
have $\Hom_S(M,S) = 0$ and $\Ext^j_S(M,S) =
\Ext^{j-1}_R(M,R)^\sigma$, as right $S$-modules, for any $j \geq
1$.}

\smallskip

Our claim is the consequence of a series of steps. First of all, for
any left $R$-modules $N$ we have the natural homomorphism of right
$S$-modules
$$
\matrix{
\Hom_R(N,R) \otimes_R S & \longrightarrow & \Hom_R(N,S) = \Hom_S(S \otimes_R
N,S)\cr
\hfill f \otimes x & \longmapsto & [n \mapsto f(n)x] \ .\hfill }
$$
It is an isomorphism if $N$ is finitely generated. Since over a
noetherian ring arbitrary direct products of (faithfully) flat
modules are (faithfully) flat and since $S$ is topologically free
over $R$ from the left as well as right it follows that $S$ is
faithfully flat as a left as well as a right $R$-module. Applying
the above isomorphism to a resolution of the $R$-module $M$ by
finitely generated projective $R$-modules therefore leads to a
natural $S$-linear isomorphism
$$
\Ext^\ast_R(M,R) \otimes_R S = \Ext^\ast_S(S \otimes_R M,S) \ .
$$
We therefore have to show the commutativity of the diagram
$$
\xymatrix{
  \Ext^j_R(M,R) \otimes_R S \ar[ddd]_{=} \ar[rr]^{\kappa^r} & &
  \Ext^j_R(M,R)^\sigma \otimes_R S \ar[d]^{=} \\
  & & \Ext^j_R({^\sigma} M,{^\sigma} R)^\sigma \otimes_R S \ar[d]^{\Ext^j_R({^\sigma} M,\sigma) \otimes id} \\
  & & \Ext^j_R({^\sigma} M,R) \otimes_R S \ar[d]^{=} \\
  \Ext^j_S(S \otimes_R M,S) \ar[rr]^{\Ext^j_S(\kappa,S)} & &
  \Ext^j_S(S \otimes_R {^\sigma} M,S)\ .   }
$$
Since all relevant modules are finitely generated this diagram is
equivalent to the diagram
$$
\xymatrix{
  \E^j_R(M,R) \otimes_R S \ar[ddd]_{=} \ar[rr]^{\kappa^r} & &
  \E^j_R(M,R)^\sigma \otimes_R S \ar[d]^{=} \\
  & & \E^j_R({^\sigma} M,{^\sigma} R)^\sigma \otimes_R S \ar[d]^{\E^j_R({^\sigma} M,\sigma) \otimes id} \\
  & & \E^j_R({^\sigma} M,R) \otimes_R S \ar[d]^{=} \\
  \E^j_S(S \otimes_R M,S) \ar[rr]^{\E^j_S(\kappa,S)} & &
  \E^j_S(S \otimes_R {^\sigma} M,S)\ .   }
$$
The vertical identifications now can be viewed as being induced
(through a projective resolution of $M$ in $\Pfr\Mfr(R)$) by the
natural maps
$$
\matrix{
\Hom_{\Pfr\Mfr(R)}(N,R) \otimes_R S & \longrightarrow &
\Hom_{\Pfr\Mfr(R)}(N,S) = \Hom_{\Pfr\Mfr(S)}(S\, \widehat{\otimes}_R\,
N,S)\cr
\hfill f \otimes x & \longmapsto & [n \mapsto f(n)x] \ .\hfill }
$$
for any $N$ in $\Pfr\Mfr(R)$. Since the whole diagram is functorial in
the pseudocompact $S$-modules $M$ we may reduce the proof of its
commutativity, by using a projective resolution of $M$ in
$\Mfr_{fg}(S)$, to the case $j = 0$, i.e., to the commutativity of the
diagram
$$
\xymatrix{
  \Hom_{\Pfr\Mfr(R)}(N,R) \otimes_R S \ar[ddd] \ar[rr]^{\kappa^r} & &
  \Hom_{\Pfr\Mfr(R)}(N,R)^\sigma \otimes_R S \ar[d]^{=} \\
  & & \Hom_{\Pfr\Mfr(R)}({^\sigma} N,{^\sigma} R)^\sigma \otimes_R S
       \ar[d]^{\Hom_{\Pfr\Mfr(R)}({^\sigma} N,\sigma) \otimes id} \\
  & & \Hom_{\Pfr\Mfr(R)}({^\sigma} N,R) \otimes_R S \ar[d] \\
  \Hom_{\Pfr\Mfr(S)}(S \otimes_R N,S) \ar[rr]^{\Hom_{\Pfr\Mfr(S)}(\kappa,S)} & &
  \Hom_{\Pfr\Mfr(S)}(S \otimes_R {^\sigma} N,S)   }
$$
for any module $N$ in $\Mfr_{fg}(S)$. This is a straightforward
explicit computation which we leave to the reader. This concludes
the proof of our claim and consequently the proof of Prop.\ 3.1.

\smallskip

To make use of Prop.\ 3.1 we need a reasonable theory of
dimension. We therefore assume now in addition that

\qquad $R$ and $S$ are (left and right) Auslander regular.

This means (for $R$) that $R$ is (left and right) regular, i.e.,
that any finitely generated (left or right) $R$-module has a
finite resolution by finitely generated projective $R$-modules.
Defining, for any finitely generated left or right $R$-module $N$,
its {\it grade} or {\it codimension} by
$$
j_R(N) := {\rm min}\{i \geq 0 : \Ext^i_R(N,R) \neq 0\}
$$
the Auslander regularity condition requires in addition that for
any such $N$ and any $k \geq 0$ we have
$$
j_R(N') \geq k \quad\hbox{for any $R$-submodule}\ N' \subseteq
\Ext^k_R(N,R)\ .
$$
A standard method to introduce a dimension filtration on the
category $\Mfr_{fg}(R)$ (cf.\ [Ste] Chap. V and [CSS]\S1) is to
define, for any $j \geq 0$, the full subcategory
$$
\matrix{
\Mfr^j(R) := & \hbox{all modules}\ M\ \hbox{in}\ \Mfr_{fg}(R)\ \hbox{such
that}\hfill\cr
 & \qquad \Ext_R^i(M',R) = 0\hfill\cr
 & \hbox{for any}\ i < j\ \hbox{and any submodule}\ M' \subseteq M\ .\hfill}
$$
It has the following properties:

-- In any short exact sequence $0 \rightarrow M' \rightarrow M
\rightarrow M'' \rightarrow 0$ in $\Mfr(R)$ the module $M$ lies in $\Mfr^j(R)$
if and only $M'$ and $M''$ lie in $\Mfr^j(R)$;

-- any module in $\Mfr_{fg}(R)$ has a unique largest submodule
contained in $\Mfr^j(R)$.

We also recall (cf.\ [CSS]\S1 and [Ve1]) that, because we are
assuming $R$ to be Auslander regular, we in fact have the
simplified characterization:

-- A module $M$ in $\Mfr_{fg}(R)$ lies in $\Mfr^j(R)$ if and only if
$\Ext_R^i(M,R) = 0$ for any $i < j$.

Note that, for any nonzero module $M$ in $\Mfr_{fg}(R)$, its grade
$j_R(M)$ is the smallest non-negative integer such that $M$ lies in
$\Mfr^{j_{R}(M)}(R)$.

Let $G_\ast(R)^{(j)}$ denote the Quillen $K$-groups of the abelian
category $\Mfr^j(R)$. In particular, $G_0(R)^{(j)}$ is the
Grothendieck group of the category $\Mfr^j(R)$ (w.r.t. short exact
sequences).

The category $\Mfr_R(S)$ can be viewed as an abelian subcategory of
both categories $\Mfr_{fg}(S)$ and $\Mfr_{fg}(R)$. Hence it becomes a
natural problem to compare the codimensions as $R$- and $S$-module,
respectively, of any $M$ in $\Mfr_R(S)$.

\smallskip

{\bf Proposition 3.2:} {\it $\Mfr^j(R) \cap \Mfr_R(S) =
\Mfr^{j+1}(S) \cap \Mfr_R(S)$ for any $j \geq 0$.}

Proof: As a consequence of Prop.\ 3.1 we have $j_S(M) = j_R(M) +
1$ for any module $M$ in $\Mfr_R(S)$.

\smallskip

The second problem which we want to address is the nature of the
natural maps
$$
G_\ast(S)^{(j+1)} \longrightarrow G_\ast(S)^{(j)}
$$
induced by the corresponding inclusion of categories. As background
for this one should keep in mind that the Gersten conjecture ([Ger])
is the claim that in the case of a commutative regular local
noetherian ring (instead of $S$) these maps are the zero maps. We
introduce the Quillen $K$-groups $G_\ast^R(S)^{(j)}$ of the category
$\Mfr^j(R) \cap \Mfr_R(S)$.

\smallskip

{\bf Lemma 3.3:} {\it For any module $N$ in $\Mfr_{fg}(R)$ we
have}
$$
j_R(N) = j_S(S \otimes_R N) \ .
$$
Proof: As discussed already after Prop.\ 3.1 we have that $S$ is
faithfully flat over $R$ and that
$$
\Ext^\ast_S(S \otimes_R N,S) = \Ext^\ast_R(N,R) \otimes_R S\ .
$$

\smallskip

{\bf Proposition 3.4:} {\it If the automorphism $\sigma$ of $R$ is
of the form $\sigma(.) = \gamma.\gamma^{-1}$ for some unit $\gamma
\in S^\times$ then the map
$$
G_\ast^R(S)^{(j)} \longrightarrow G_\ast(S)^{(j)}
$$
induced by the inclusion $\Mfr^j(R) \cap \Mfr_R(S) \subseteq
\Mfr^j(S)$ is the zero map for any $j \geq 0$.}

Proof: Suppose that the module $M$ lies in $\Mfr^j(R) \cap
\Mfr_R(S)$. Then by Prop.\ 3.2 and Lemma 3.3 all three $S$-modules
which constitute the exact sequence of Remark 2.3.i lie in the
abelian category $\Mfr^j(S)$. We therefore may view this short
exact sequence as an exact sequence of natural transformations
between the exact functors from the category $\Mfr^j(R) \cap
\Mfr_R(S)$ into the category $\Mfr^j(S)$ given by inclusion and $S
\otimes_R .$, respectively. In this situation a fundamental
theorem of Quillen asserts ([Qui] \S3 Cor. 1) that the map induced
on $K$-groups by the third functor is the difference of the maps
induced by the first two functors. But these first two functors
coincide. Hence the inclusion functor induces the zero map.

\smallskip

Since by Prop.\ 3.2 we have the commutative triangle
$$
\xymatrix{
                & G_\ast^R(S)^{(j)} \ar[dl] \ar[dr]^{0}             \\
 G_\ast(S)^{(j+1)}  \ar[rr] & &    G_\ast(S)^{(j)} }
$$
Prop.\ 3.4 may be viewed as a partial answer to our problem. But
it is not clear whether this vanishing result can be considered as
evidence for some version of Gersten's conjecture for certain
noncommutative rings.

\medskip

{\bf 4. Application to Iwasawa algebras}

\medskip

Let $G$ be a compact $p$-adic Lie group. Its so called Iwasawa
algebra
$$
 \Lambda(G) := \mathop{\rm lim}\limits_{\longleftarrow}
 \Zdss_p[G/U]
$$
is the projective limit of the group rings $\Zdss_p[G/U]$ where
$U$ ranges over the open normal subgroups in $G$. It is known to
be (left and right) noetherian ([Laz] (V.2.2.4)). An important
feature of $\Lambda(G)$ is that it carries a natural compact
Hausdorff topology. A fundamental system of open neighbourhoods of
zero in $\Lambda(G)$ is given by the submodules $p^n\Lambda(G) +
I(U)$ with $n \in \Ndss$ and $U$ as before, where $I(U)$ denotes
the kernel of the projection map $\Lambda(G) \rightarrow
\Lambda(G/U)$ (cf.\ [NSW] Prop.\ (5.2.17) or [Laz] (II.2.2.2)). It
follows that $\Lambda(G)$ is left and right pseudocompact.

We remark that if $G$ is pro-$p$ then $\Lambda(G)$ is local, and if in
addition $G$ has no element of order $p$ then $\Lambda(G)$ is a
regular local integral domain ([Neu]).

We now assume that there is a closed normal subgroup $H \subseteq
G$ such that $G/H \cong \Zdss_p$. Then with the pair of
pseudocompact rings $R := \Lambda(H) \subseteq \Lambda(G) =: S$ we
are in the situation considered in this paper. This is seen by
noting that $G$ necessarily is the semidirect product of $H$ and a
subgroup $\Gamma \cong \Zdss_p$. We pick once and for all a
topological generator $\gamma$ of $\Gamma$ and define the element
$t := \gamma - 1$ in $\Lambda(\Gamma) \subseteq \Lambda(G)$. It is
well known that $\Lambda(\Gamma)$ is the formal power series ring
in the variable $t$ over $\Zdss_p$. Hence (use Lemma 1.6) the
Iwasawa algebra $S = \Lambda(G)$ is the skew power series ring
$\Lambda(H)[[t;\sigma,\delta]]$ where $\sigma(.) := \gamma
.\gamma^{-1}$ and $\delta := \sigma - id$ (cf.\ [Ve3] Example
5.1).

We recall (cf.\ [CSS]\S1 and [Ve1]) that $G$ always contains an
open subgroup $G'$ such that $\Lambda(G')$ is left and right
Auslander regular. Since
$$
\Ext_{\Lambda(G)}^i(M,\Lambda(G)) = \Ext_{\Lambda(G')}^i(M,\Lambda(G'))
$$
for any $i \geq 0$ and any module $M$ in $\Mfr_{fg}(\Lambda(G))$
(cf.\ [Jan] Lemma 2.3) the dimension theory of section 3 continues
to hold for general $\Lambda(G)$ even if it is not Auslander
regular.

Hence our main results Prop.s 2.2, 3.1, 3.2, and 3.4 hold true for
the pair $\Lambda(H) \subseteq \Lambda(G)$. In Iwasawa theory the
modules in $\Mfr^2(\Lambda(G))$ are called pseudo-null modules. We
conclude by pointing out that, as a special case of our results,
any pseudo-null $\Lambda(G)$-module which is finitely generated as
a $\Lambda(H)$-module has zero class in the Grothendieck group
$G_0(\Lambda(G))^{(1)}$ of all finitely generated torsion
$\Lambda(G)$-modules.

\medskip

{\bf 5. Skew Laurent series rings}

\medskip

Let $R$ be again a pseudocompact ring and $\sigma$ a topological
automorphism of $R$. The left $R$-module $T$ of all (left) formal
Laurent series $\sum_{i \gg -\infty} a_i t^i$ over $R$ carries a
unique unital ring structure which extends the $R$-module
structure and satisfies $ta=\sigma(a)t$. The ring
$R[[t,t^{-1};\sigma]] := T$ is called a skew Laurent series ring.
It can alternatively be described by right Laurent series, it
contains $S := R[[t;\sigma]]$ as a subring, and it also arises as
the localization $S_t$ of $S$ with respect to the Ore set
$\{1,t,t^2,\ldots\}$. In particular, $T$ is a flat (left and
right) $S$-module. In the following we assume that $R$ is
noetherian. Then $S$, by [MCR] 1.4.5, and, as a localization of a
noetherian ring, $T$ both are noetherian as well.

Similarly as before we let $\Mfr_R(T)$ denote the abelian category
of left $T$-modules which are finitely generated as $R$-modules.
Any such module $M$ lies in $\Mfr_R(S)$ and satisfies $T \otimes_S
M = M$. From this observation it is immediate that, for any $M$ in
$\Mfr_R(T)$, we obtain an analogue of Prop.\ 2.2 in the form of an
exact sequence of (left) $T$-modules
$$
0 \longrightarrow T \otimes_R M \longrightarrow T \otimes_R M
\longrightarrow M \longrightarrow 0\ .
$$
Indeed, consider $M$ as an $S$-module, apply $T\otimes_S .$ to the
exact sequence resulting from that proposition and observe that an
analog of Remark 2.3.i applies because $\sigma(.)=t . t^{-1}$.

Furthermore, for $M$ in $\Mfr_R(T)$, multiplication by $t$ on the
right $S$-module $\Ext^j_R(M,R)^\sigma$ is easily checked to be
bijective. Hence the $S$-module structure extends uniquely to a
right $T$-module structure on $\Ext^j_R(M,R)^\sigma$. Using the
natural isomorphism $\Ext^j_T(M,T)\cong \Ext^j_S(M,S)\otimes_S T$
of right $T$-modules we obtain, again by base change, the analogue
of Prop.\ 3.1:
$$
\Ext^j_T(M,T)\cong \Ext^{j-1}_R(M,R)^\sigma
$$
as right $T$-modules.

Now assume that $R$ in addition is Auslander regular. By [LVO]
III.3.4.6(3) the ring $S$ then also is Auslander regular. The
increasing filtration $(St^{-j})_{j \in \Zdss}$ on $T$ is complete
with associated graded ring isomorphic to $R[t,t^{-1};\sigma]$. It
therefore follows from [LVO] III.3.4.6(2) and III.2.2.5 that $T$
is Auslander regular. We also note that $T$ is faithfully flat
over $R$. Thus the analogs of Prop.\ 3.2 and Prop.\ 3.4 (the
assumption there is satisfied in the present situation) hold for
the ring extension $T$ over $R$ by exactly the same arguments.

Finally, we point out that all the results of this section also
hold for skew Laurent polynomial rings $R[t,t^{-1};\sigma]$ over
arbitrary (noetherian, respectively Auslander regular) rings $R$
with ring automorphism $\sigma$. We omit the proofs since they are
analogous but simpler.

\bigskip

{\bf References}

\parindent=23truept

\ref{[Bru]} Brumer A.: Pseudocompact Algebras, Profinite Groups and
Class Formations. J. Algebra 4, 442-470 (1966)

\ref{[CH]} Chatters A.W., Hajarnavis C.R.: Rings with chain conditions.
Boston-London-Melbourne: Pitman 1980

\ref{[CSS]} Coates J., Schneider P., Sujatha R.:
Modules over Iwasawa algebras. J. Inst. Math. Jussieu 2, 73--108
(2003)


\ref{[Ga1]} Gabriel P.: Des Cat\'egories Ab\'eliennes. Bull. Soc.
math. France 90, 323-448 (1962)

\ref{[Ga2]} Gabriel P.: Etude infinitesimale des schemas en
groupes. In Sch\'emas en Groupe I (SGA3), Exp. VII B. Lecture
Notes Math. 151, pp. 4746-562. Berlin-Heidelberg-New York:
Springer 1970

\ref{[Ger]} Gersten S.M.: Some exact sequences in the higher
$K$-theory of rings. In Algebraic $K$-theory I. Lecture Notes Math.
341, pp. 211-243. Berlin-Heidelberg-New York: Springer 1973

\ref{[Jan]} Jannsen U.: Iwasawa modules up to isomorphism. Adv.
Studies Pure Math. 17, 171-207 (1989)

\ref{[Laz]} Lazard M.: Groupes analytiques $p$-adique. Publ. Math.
IHES 26, 389-603 (1965)

\ref{[LVO]} Li Huishi, van Oystaeyen F.: Zariskian Filtrations.
Dordrecht: Kluwer 1996

\ref{[MCR]} McConnell J.C., Robson J.C.: Noncommutative
noetherian rings. Chi- chester: Wiley 1987

\ref{[NSW]} Neukirch J., Schmidt A., Wingberg K.: Cohomology of Number
Fields. Grundlehren math. Wiss. 323.  Berlin-Heidelberg-New York:
Springer 2000

\ref{[Neu]} Neumann A.: Completed group algebras without zero
divisors. Arch. Math. 51, 496-499 (1988)


\ref{[Qui]} Quillen D.: Higher algebraic $K$-theory I. In Algebraic
$K$-theory I. Lecture Notes Math. 341, pp. 85-139.
Berlin-Heidelberg-New York: Springer 1973

\ref{[Ree]} Rees D.: A theorem of homological algebra. Proc. Cambridge
Philos. Soc. 52, 605-610 (1956)

\ref{[Ste]} Stenstr{\"o}m B.: Rings of Quotients.
Berlin-Heidelberg-New York: Sprin- ger 1975

\ref{[vGB]} van Gastel M., van den Bergh M.: Graded Modules of
Gelfand-Kirillov Dimension One over Three-Dimensional Artin-Shelter
Regular Algebras. J. Algebra 196, 251-282 (1997)

\ref{[Ve1]} Venjakob O.: On the structure theory of the Iwasawa algebra of a
$p$-adic Lie group. J. Eur. Math. Soc. 4, 271-311 (2002)

\ref{[Ve2]} Venjakob O.: A non-commutative Weierstrass
preparation theorem and applications to Iwasawa theory. J. reine
angew. Math. 559, 153-191 (2003)

\ref{[Ve3]} Venjakob O.: Characteristic Elements in Noncommutative
Iwasawa Theory. Habilitationsschrift, Heidelberg 2003

\parindent=0pt

\bigskip

Universit\"{a}t M{\"u}nster, Mathematisches Institut, Einsteinstr. 62,
48291 M{\"u}nster, Germany\hfill\break
 pschnei@math.uni-muenster.de\hfill\break
 http://www.uni-muenster.de/math/u/schneider/

Universit\"{a}t Heidelberg, Mathematisches Institut, Im Neuenheimer
Feld 288, 69120 Heidelberg, Germany\hfill\break
 otmar@mathi.uni-heidelberg.de\hfill\break
 http://www.mathi.uni-heidelberg.de/$\,\tilde{}\,$otmar/

\end